\documentclass[12pt]{article}
\usepackage{amssymb,amsthm,amsmath,latexsym}
\usepackage{url}
\usepackage{color}
\usepackage{comment}
\usepackage{lscape}
\usepackage{mathtools}
\usepackage{enumitem}
\usepackage{graphicx}
\newtheorem{thm}{Theorem}[section]
\newtheorem{prop}[thm]{Proposition}
\newtheorem{lem}[thm]{Lemma}
\newtheorem{cor}[thm]{Corollary}
\theoremstyle{remark}
\newtheorem{rem}[thm]{Remark}
\newcommand{\FF}{\mathbb{F}}
\newcommand{\ZZ}{\mathbb{Z}}
\newcommand{\RR}{\mathbb{R}}
\newcommand{\0}{\mathbf{0}}

\allowdisplaybreaks

\begin{document}

\title{Optimal odd unimodular lattices in dimensions $49$ through $59$}

\author{
Masaaki Harada\thanks{
Research Center for Pure and Applied Mathematics,
Graduate School of Information Sciences,
Tohoku University, Sendai 980--8579, Japan.
email: \protect\url{mharada@tohoku.ac.jp}
\protect\newline
\protect\newline
}
}

\maketitle

\begin{abstract}
It is shown that the largest minimum norm among unimodular lattices
is $4$ in dimension $49$ and $5$ in dimensions $50$, $51$, $57$ and $58$.
We also study the theta series of optimal unimodular lattices
in every dimension from $49$ through $59$.
In dimensions $50$, $51$ and $57$, the first known examples of
unimodular lattices with minimum norm $5$ are obtained 
from self-dual codes by Construction~A.
There exist at least $25$, $2$, $3$, $7$ and $12$ nonisomorphic
optimal unimodular lattices in dimensions $49$, $55$, $57$, $58$ and $59$,
respectively.
\end{abstract}

\section{Introduction}\label{Sec:1}

A (Euclidean) integral lattice $L$ in dimension $n$
is \emph{unimodular} if $L = L^{*}$, where
$L^{*}$ denotes the dual lattice of $L$.
Unimodular lattices are studied from several
viewpoints (see~\cite{SPLAG} for an extensive bibliography).
Self-dual codes and unimodular lattices are closely related through
Construction~A.

It is a fundamental problem to determine the largest minimum norm $\mu(n)$
among all unimodular lattices in dimension $n$.
A unimodular lattice with the largest minimum norm in that dimension is called 
\emph{optimal}~\cite{CS98}.
Bounds on the largest minimum norms in dimensions $n \le 40$
are given in~\cite[Table 1]{CS98}.
This work was extended in~\cite{NG}.
Recently, Kubischta~\cite{Kubischta} constructed a unimodular
lattice with minimum norm $4$ in dimension $37$.
This completes the determination of the largest minimum norms $\mu(n)$ 
for $n \le 48$.
To extend this determination through dimension $61$, it suffices to
settle five cases in view of the known bounds and constructions
(see Table~\ref{Tab:mu}).
In dimensions $49$ and $58$, the lattices constructed
in~\cite{SYS} give the lower bounds $\mu(49)\geq4$ and
$\mu(58)\geq5$, respectively. 
We establish the matching upper bounds.
In dimensions $50$, $51$ and $57$, we construct the first known
examples of unimodular lattices attaining the upper bound 
given in~\cite[Table~2]{Ga04}.
Hence, we conclude that $\mu(49)=4$ and $\mu(n)=5$ for $n\in\{50,51,57,58\}$.

This paper is organized as follows.
In Section~\ref{Sec:2}, we review the definitions, notation, and basic
results used throughout the paper.
In Section~\ref{Sec:T}, we prove the upper bounds $\mu(49)\leq4$ and
$\mu(58)\leq5$. We also study the theta series of optimal unimodular lattices
in every dimension from $49$ through $59$.
In Section~\ref{Sec:C}, we apply Construction~A to self-dual codes
to obtain optimal unimodular lattices in dimensions
$49$, $50$, $51$, $55$, $57$, $58$ and $59$.
The lattice constructed in dimension $50$ is also $s$-extremal.
Finally, we show that there exist at least $25$, $2$, $3$, $7$ and $12$
nonisomorphic optimal unimodular lattices in dimensions
$49$, $55$, $57$, $58$ and $59$, respectively.

\section{Preliminaries}\label{Sec:2}

In this section, we recall the definitions, notation and basic results
used throughout the paper.

\subsection{Unimodular lattices}
A (Euclidean) lattice $L \subset \RR^n$ 
in dimension $n$
is \emph{unimodular} if
$L = L^{*}$, where
the dual lattice $L^{*}$ of $L$ is defined as
$\{ x \in {\RR}^n \mid (x,y) \in \ZZ \text{ for all }
y \in L\}$ under the standard inner product $(x,y)$.
The \emph{Gram matrix} of vectors $v_1,v_2,\ldots,v_m\in\RR^n$ is the
matrix $\bigl((v_i,v_j)\bigr)_{1\leq i,j\leq m}$.
If these vectors form a $\ZZ$-basis of a lattice $L$,
this matrix is called a Gram matrix of $L$.
A unimodular lattice $L$ is \emph{even} 
if the norm $(x,x)$ of every vector $x$ of $L$ is even,
and \emph{odd} otherwise.
An even unimodular lattice in dimension $n$
exists if and only if $n \equiv 0 \pmod 8$, while
an odd unimodular lattice exists for every dimension.
The \emph{minimum norm} $\min(L)$ of 
a unimodular lattice $L$ is the smallest 
norm among all nonzero vectors of $L$.
The \emph{kissing number} of $L$ is the number of vectors of minimum norm in $L$.
Two lattices $L$ and $L'$ are \emph{isomorphic}
if there is an orthogonal matrix $A$ with
$L' = L \cdot A$, where
$L \cdot A =\{xA \mid x \in L\}$.

Let $L$ be an odd unimodular lattice in dimension $n$.
The \emph{shadow} $S(L)$ of $L$ is defined to be $S(L)= L_0^* \setminus L$,
where $L_0$ denotes the even sublattice of $L$.
Since $|L:L_0|=|L_0^*:L|=2$, there are cosets
$\mathcal C_1,\mathcal C_2,\mathcal C_3$ of $L_0$ such that
\begin{equation*}
\begin{aligned}
L_0^*&=L_0\sqcup \mathcal C_1\sqcup \mathcal C_2\sqcup \mathcal C_3,\\
L&=L_0\sqcup \mathcal C_2 \text{ and}\\
S(L)&=\mathcal C_1\sqcup \mathcal C_3,
\end{aligned}
\end{equation*}
where $\sqcup$ denotes a disjoint union.
Since the quotient group $L_0^*/L$ has order $2$, we see that
$2s\in L$ for every $s\in S(L)$.
Shadows of odd unimodular lattices appeared in~\cite{CS98} 
(see also~\cite[p.~440]{SPLAG}),
and shadows play an important role in the study of odd unimodular lattices.

The \emph{theta series} of an odd unimodular lattice $L$ and its shadow $S(L)$
are the formal power series
$\theta_{L}(q) = \sum_{x \in L} q^{(x,x)}$ and
$\theta_{S(L)}(q) = \sum_{x \in S(L)} q^{(x,x)}$, respectively.
Conway and Sloane~\cite{CS98} showed that
when the theta series of an odd unimodular lattice $L$
in dimension $n$
is written as
\begin{equation}
\label{eq:theta}
 \sum_{j =0}^{\lfloor n/8\rfloor} a_j\theta_3(q)^{n-8j}\Delta_8(q)^j,
\end{equation}
the theta series of the shadow $S(L)$
is written as
\begin{equation}
\label{eq:theta-S}
\sum_{j=0}^{\lfloor n/8\rfloor}
\frac{(-1)^j}{16^j} a_j\theta_2(q)^{n-8j}\theta_4(q^2)^{8j},
\end{equation}
where
$\Delta_8(q) = q \prod_{m=1}^{\infty} (1 - q^{2m-1})^8(1-q^{4m})^8$,
and $\theta_2(q)$, $\theta_3(q)$ and $\theta_4(q)$ are the Jacobi
theta series~\cite{SPLAG}.
By considering the theta series of unimodular lattices and their shadows,
Rains and Sloane~\cite{RS} showed that a unimodular lattice $L$
in dimension $n$ satisfies
\begin{equation}\label{eq:RS}
\min(L)\le 2\left\lfloor\frac{n}{24}\right\rfloor+2,
\end{equation}
unless $n=23$, in which case $\min(L)\leq3$.

\subsection{Self-dual $\ZZ_k$-codes and Construction~A}
Let $\ZZ_{k}$ be the ring 
of integers modulo $k$, where $k$ 
is a positive integer. 
In this paper, we always assume that $k\geq 2$ and 
we take the set $\ZZ_{k}$ to be $\{0,1,\ldots,k-1\}$.
A $\ZZ_{k}$-code $C$ of length $n$
is a $\ZZ_{k}$-submodule of $\ZZ_{k}^n$.
A $\ZZ_{k}$-code $C$ is \emph{self-dual} if $C=C^\perp$, where
the dual code $C^\perp$ of $C$ is defined as 
\[
\{ x \in \ZZ_{k}^n \mid x \cdot y = 0\text{ for all }y \in C\}
\]
under the standard inner product $x \cdot y$. 

Construction~A is a standard method for constructing unimodular
lattices from self-dual codes.
Many unimodular lattices with large minimum norms 
have been obtained by this construction.
It allows algebraic constructions and computational searches
for self-dual codes to be used in the search for optimal
unimodular lattices.
In particular, if $C$ is a self-dual $\ZZ_{k}$-code of length $n$, then 
the following lattice 
\[
A_{k}(C)=
\frac{1}{\sqrt{k}}\{x \in \ZZ^n \mid x \pmod k \in C\}
\]
is a unimodular lattice in dimension $n$,
where reduction modulo $k$ is taken coordinatewise.

\section{Upper bounds on minimum norms and theta series}\label{Sec:T}

In this section, it is shown that $\mu(49)=4$ and $\mu(58)=5$.
We also study the theta series of optimal unimodular lattices
in every dimension from $49$ through $59$.
We used \textsc{Magma}~\cite{Magma} to verify the coefficients
in the theta series expansions displayed in this section.

\subsection{Preliminary results}

Let $L$ be an odd unimodular lattice.
Let $L_0$ and $S(L)$ denote its even sublattice and its shadow,
respectively.
For $r>0$, let
\[
 L_r=\{x\in L\mid(x,x)=r\}\text{ and }
 S(L)_r=\{s\in S(L)\mid(s,s)=r\}.
\]
Each of these sets is a disjoint union of pairs $\{x,-x\}$.
In particular, $|S(L)_r|/2$ is a nonnegative integer.

\begin{lem}\label{lem:theta-common-coefficients}
Suppose that $n\in\{49,50,51,57,58,59\}$.
Let $L$ be a unimodular lattice in dimension $n$.
Suppose that $\min(L)\geq4$ if $n=49$ and $\min(L)\geq5$ otherwise.
Let
\[
 \ell=\left\lfloor\frac n8\right\rfloor,
 \rho=\frac{n-8\ell}{4}+2,
 \beta=\frac{|S(L)_\rho|}{2}\text{ and }
 \gamma=\frac{|S(L)_{\rho+2}|}{2}.
\]
Then $S(L)_{\rho-2}=\emptyset$, $a_0=1$, and the coefficients
$a_1,a_2,\ldots,a_\ell$ in~\eqref{eq:theta} are given in
Table~\ref{Tab:theta-coefficients}.
\end{lem}

\begin{proof}
The constant term in~\eqref{eq:theta} gives $a_0=1$.
Using
\begin{align*}
 \theta_3(q)&=1+2q+2q^4+\cdots\text{ and }\\
 \Delta_8(q)&=q-8q^2+28q^3-64q^4+126q^5+\cdots,
\end{align*}
we set the coefficient of $q^r$ equal to zero for each integer $r$
with $1\leq r<4$ if $n=49$, and with $1\leq r<5$ otherwise.
This determines $a_1,a_2,a_3$ for $n=49$ and $a_1,a_2,a_3,a_4$
for the other five dimensions, as listed in Table~\ref{Tab:theta-coefficients}.

\begin{table}[thb]
\caption{Coefficients in~\eqref{eq:theta}}
\label{Tab:theta-coefficients}
\centering
\medskip
{\small
\begin{tabular}{c|c|c|c|c|c|c|c}
\noalign{\hrule height1pt}
$n$&{$a_1$}&{$a_2$}&{$a_3$}&{$a_4$}&{$a_5$}&{$a_6$}&{$a_7$}\\
\hline
$49$&$-98$&$2548$&$-14896$&$71\beta+\gamma$&$-4096\beta$&$0$&--\\
$50$&$-100$&$2700$&$-17200$&$9100$&$-2048\beta$&$0$&--\\
$51$&$-102$&$2856$&$-19720$&$13770$&$-1024\beta$&$0$&--\\
$57$&$-114$&$3876$&$-39824$&$77406$&$-16(87\beta+\gamma)$&$65536\beta$&$0$\\
$58$&$-116$&$4060$&$-44080$&$96860$&$-688\beta-8\gamma$&$32768\beta$&$0$\\
$59$&$-118$&$4248$&$-48616$&$119770$&$-340\beta-4\gamma$&$16384\beta$&$0$\\
\noalign{\hrule height1pt}
\end{tabular}
}
\end{table}

We have $4(\rho-2)=1,2,3,1,2,3$ for $n=49,50,51,57,58,59$, respectively.
If $s\in S(L)_{\rho-2}$, then $2s$ is a nonzero vector of $L$ of norm
$4(\rho-2)$, which is smaller than the assumed minimum norm.
Hence, $S(L)_{\rho-2}=\emptyset$.
Only the term with $j=\ell$ can contribute to the coefficient of
$q^{\rho-2}$ in~\eqref{eq:theta-S}, and its leading factor is nonzero.
Hence, we conclude that $a_\ell=0$.
Using 
\[
\theta_2(q)=2q^{1/4}(1+q^2+\cdots)\text{ and }
\theta_4(q^2)=1-2q^2+\cdots, 
\]
we compare the coefficients of
$q^\rho$ and $q^{\rho+2}$ and we obtain
\begin{align*}
 2\beta&=(-1)^{\ell-1}2^{n-12\ell+12}a_{\ell-1}\text{ and }\\
 2\gamma&=(-1)^{\ell-2}2^{n-12\ell+24}a_{\ell-2}
          -(24\ell-n-24)\cdot2\beta.
\end{align*}
For $n=49$, these equations give $a_4$ and $a_5$.
For $n=57,58,59$, they give $a_5$ and $a_6$.
For $n=50$ and $51$, the first equation gives $a_5$.
Since $a_4$ is already known in these two dimensions,
the second equation determines $\gamma$ in terms of $\beta$.
\end{proof}

\begin{rem}
It follows from Table~\ref{Tab:theta-coefficients} and~\eqref{eq:theta}
that the entire theta series $\theta_L(q)$ is determined by $\beta$
for $n=50,51$ and by $\beta$ and $\gamma$ for $n=49,57,58,59$.
\end{rem}

The following lemma gives the restrictions obtained from inner products
between short vectors of the shadow. 
The argument in~(i) is similar to that in~\cite[p.~358]{CS98}.
Let $J_n$ denote the all-one matrix of order $n$.

\begin{lem}\label{lem:theta-shadow-bounds}
Let $L$ be an odd unimodular lattice in dimension $n$.
\begin{enumerate}
\item\label{item:shadow-five-halves}
If $\min(L)\geq5$, then $|S(L)_{5/2}|/2\in\{0,1,2\}$.
\item\label{item:shadow-eleven-quarters}
If $\min(L)\geq5$, then $|S(L)_{11/4}|/2\in\{0,1,\ldots,10,12\}$.
\item\label{item:shadow-nine-quarters}
If $\min(L)\geq4$, then $|S(L)_{9/4}|/2\leq n$.
Let $\beta=|S(L)_{9/4}|/2$.  Then the set
$S(L)_{9/4}\cap \mathcal C_1$ has $\beta$ elements and Gram matrix
$2I_\beta+J_\beta/4$.
\end{enumerate}
\end{lem}

\begin{proof}
For the assertion~(i), suppose that $s,t\in S(L)_{5/2}$ belong to the same
coset of $L_0$ and $t\ne\pm s$.
Then we have $s-t\in L_0$. Since $2t\in L$ has norm $10$, we have
$2t\in L_0$ and hence we have $s+t\in L_0$.
Since both vectors are nonzero, the evenness of $L_0$ and $\min(L)\geq5$ give
\[
 12\leq(s-t,s-t)+(s+t,s+t)=10,
\]
a contradiction.
Thus, each of the two shadow cosets contains at most one pair
$\{s,-s\}$ of norm $5/2$, which proves the assertion~(i).

Set $r=11/4$ in~(ii) and $r=9/4$ in~(iii).
Let $\beta=|S(L)_r|/2$ and let
\[
  X=S(L)_r\cap\mathcal C_1.
\]
For $s\in S(L)_r$, we have $2s\in L$ and $(2s,2s)=4r$.
Since $4r=11$ in~(ii) and $4r=9$ in~(iii), we see that $2s\notin L_0$.
Thus, $s$ and $-s$ lie in different cosets of $L_0$,
since $s-(-s)=2s$.
As $\mathcal C_1$ and $\mathcal C_3$ are the two shadow cosets,
exactly one vector from each pair $\{s,-s\}$ belongs to $X$.
Hence, we have $|X|=|S(L)_r|/2=\beta$.
For distinct $s,t\in X$, we have
\[
 s-t\in L_0\text{ and }s+t\in L\setminus L_0.
\]
Since both vectors are nonzero, their norms are at least $6$ and $5$
when $r=11/4$, and at least $4$ and $5$ when $r=9/4$.
Since the sum of their norms is $4r$, we obtain
\[
 (s,t)=
 \begin{cases}
 -1/4&\text{if }r=11/4,\\
  1/4&\text{if }r=9/4.
 \end{cases}
\]

For $r=11/4$, we have
\[
 0\leq\left(\sum_{s\in X}s,\sum_{s\in X}s\right)
 =\frac{11\beta}{4}-\frac{\beta(\beta-1)}4
 =\frac{\beta(12-\beta)}4.
\]
This gives $\beta\leq12$.
Suppose that $\beta=11$ and let $v=-\sum_{s\in X}s$.
Then we obtain $(v,v)=11/4$.
Choose a vector $s_0\in X$. Since $2s_0\in L$ and $(4s_0,4s_0)=44$,
we see that $4s_0\in L_0$. Therefore, we obtain
\[
 v-s_0=-12s_0-\sum_{s\in X}(s-s_0)\in L_0.
\]
This gives $v\in X$.
On the other hand, for every vector $t\in X$, we have
\[
 (v,t)=-\left(\frac{11}{4}+10\left(-\frac14\right)\right)=-\frac14.
\]
Since $v\in X$, taking $t=v$ gives $(v,v)=-1/4$,
contrary to $(v,v)=11/4$.
Hence, $\beta\neq11$, which proves the assertion~(ii).

For $r=9/4$ and $\beta>0$, the Gram matrix of $X$ is
$G=2I_\beta+J_\beta/4$.
For every nonzero vector $v=(v_1,v_2,\ldots,v_\beta)\in\RR^\beta$, we have
\[
 vGv^T=2\sum_{i=1}^{\beta}v_i^2
       +\frac14\left(\sum_{i=1}^{\beta}v_i\right)^2>0.
\]
Thus, the vectors in $X$ are linearly independent.  This gives $\beta\leq n$.
The assertion is immediate when $\beta=0$, which completes the proof.
\end{proof}

Let $P(x)$ be a homogeneous harmonic polynomial of degree $2$ on
$\RR^n$ (see~\cite[Section~1]{Pache05}).
Let $L$ be an odd unimodular lattice in dimension $n$.
For $X=L$ or $S(L)$, define
\[
 \theta_{X,P}(q)=\sum_{x\in X}P(x)q^{(x,x)}.
\]
Let 
\[
\Phi(q)=\theta_4(q)^4-\theta_2(q)^4.
\]
By Proposition~8 and Theorem~12 in~\cite{Pache05},
$\theta_{L,P}(q)$ is a modular form of signature $(2,n/2+2,-)$
and can be expressed as $\Phi$ times a polynomial in
$\theta_3$ and $\Delta_8$.
Since their weights are $2$, $1/2$ and $4$, respectively,
a monomial $\Phi\theta_3^r\Delta_8^j$ occurring in this expression
must satisfy $2+r/2+4j=n/2+2$.
Hence, there exist coefficients $c_j(P)$ such that
\begin{equation}\label{eq:common-weighted-L}
 \theta_{L,P}(q)
 =\Phi(q)\sum_{j=0}^{\lfloor n/8\rfloor}
   c_j(P)\theta_3(q)^{n-8j}\Delta_8(q)^j.
\end{equation}
By Proposition~10 in~\cite{Pache05} and the transformation formulas
in~\cite[p.~1096]{Pache05}, we obtain
\begin{equation}\label{eq:common-weighted-S}
 \theta_{S(L),P}(q)
 =-\bigl(\theta_3(q)^4+\theta_4(q)^4\bigr)
  \sum_{j=0}^{\lfloor n/8\rfloor}
  \frac{(-1)^j c_j(P)}{16^j}
  \theta_2(q)^{n-8j}\theta_4(q^2)^{8j}.
\end{equation}
We use the following identity in dimensions $49$, $58$ and $59$.

\begin{lem}\label{lem:theta-weighted-identities}
Let $L$ be an odd unimodular lattice in dimension $n$.
Suppose that $\min(L)\geq m$, where the parameters 
$(m,\rho,k,\lambda)$ are given by
\[
(m,\rho,k,\lambda)=
\begin{cases}
(4,9/4,47,-4)  &\text{if }n=49,\\
(5,5/2,62,1/2) &\text{if }n=58,\\
(5,11/4,61,1)   &\text{if }n=59.
\end{cases}
\]
For every homogeneous harmonic polynomial $P$ of degree $2$ on $\RR^n$,
\begin{equation}\label{eq:theta-common-weighted-relation}
 \sum_{s\in S(L)_{\rho+2}}P(s)+k\sum_{s\in S(L)_\rho}P(s)
 =\lambda\sum_{y\in L_m}P(y).
\end{equation}
\end{lem}

\begin{proof}
In~\eqref{eq:common-weighted-L}, the factor multiplying $c_j(P)$ has
leading term $q^j$. Since $P(0)=0$ and $\min(L)\geq m$, comparing
the coefficients of $1,q,\ldots,q^{m-1}$ gives
$c_0(P)=\cdots=c_{m-1}(P)=0$. Comparing the coefficients of $q^m$ gives
\[
 c_m(P)=\sum_{y\in L_m}P(y).
\]

Using
\begin{align*}
 \theta_3(q)^4+\theta_4(q)^4&=2(1+24q^2+\cdots),\\
 \theta_2(q)&=2q^{1/4}(1+q^2+\cdots)\text{ and}\\
 \theta_4(q^2)&=1-2q^2+\cdots,
\end{align*}
we see that the term with index $j$ in~\eqref{eq:common-weighted-S} is
\[
 (-1)^{j+1}2^{n-12j+1}c_j(P)\,q^{(n-8j)/4}
 \bigl(1+(24+n-24j)q^2+\cdots\bigr).
\]
Here the coefficient $24+n-24j$ is the sum of $24$, $n-8j$ and $-16j$
from the three factors, respectively.

For each of the three dimensions $n=49,58,59$, we have
$\lfloor n/8\rfloor=m+2$ and $\rho=(n-8(m+1))/4$.
Thus, only the term with $j=m+2$ can contribute to the coefficient
of $q^{\rho-2}$. Since $S(L)_{\rho-2}=\emptyset$, by
Lemma~\ref{lem:theta-common-coefficients}, this gives $c_{m+2}(P)=0$.
It remains to consider the terms with $j=m$ and $j=m+1$.
The term with $j=m$ starts with $\lambda c_m(P)q^{\rho+2}$,
where $\lambda=(-1)^{m+1}2^{n-12m+1}$.
Let $d(P)$ denote the coefficient of $q^\rho$ in the term with $j=m+1$.
Since $24+n-24(m+1)=n-24m=-k$, we obtain
\[
 \theta_{S(L),P}(q)
 =d(P)q^\rho+\bigl(\lambda c_m(P)-k d(P)\bigr)q^{\rho+2}+\cdots.
\]
The coefficients of $q^\rho$ and $q^{\rho+2}$ are the sums of $P(s)$
over $S(L)_\rho$ and $S(L)_{\rho+2}$, respectively.
Adding $k$ times the first coefficient to the second eliminates $d(P)$
and gives~\eqref{eq:theta-common-weighted-relation}.
\end{proof}

We also use the following congruence in these three dimensions.
Let $L$ be an odd unimodular lattice and let $u\in L$ have odd norm.
Then we have $(s,u)\in\ZZ+1/2$ for every vector $s\in S(L)$.
Indeed, $2u\in L_0$ gives $2(s,u)\in\ZZ$, while
$(s,u)\in\ZZ$ implies $s\in L^*=L$, since
$L=L_0\cup(u+L_0)$.
Let $(s,u)=a+1/2$ with $a\in\ZZ$.  Then we have
\[
 (s,u)^2-\frac14=a(a+1)\in2\ZZ.
\]
Therefore, for each pair $\{s,-s\}$, its contribution is 
\[
 \left((s,u)^2-\frac14\right)
 +\left((-s,u)^2-\frac14\right)\in4\ZZ.
\]
Since $S(L)_r$ is a disjoint union of pairs $\{s,-s\}$ for $r>0$, we conclude that
\begin{equation}\label{eq:theta-shadow-square-congruence}
 \sum_{s\in S(L)_r}\left((s,u)^2-\frac14\right)\equiv0\pmod4.
\end{equation}

\subsection{Dimension 49}

Recall that $\mu(n)$ denotes the largest minimum norm 
among all unimodular lattices in dimension $n$.
The upper bound $\mu(49)\leq5$ is given in~\cite[Table~2]{Ga04}.

\begin{prop}\label{prop:rank49-min4-theta}
Let $L$ be a unimodular lattice in dimension $49$.
If $\min(L)\geq4$, then its theta series has the form
\begin{equation}\label{eq:rank49-min4-theta}
\begin{aligned}
 \theta_L(q)
 ={}&1+(71\beta+\gamma-5390)q^4
     +(305760-3954\beta+2\gamma)q^5\\
 &+(23582720+86704\beta-48\gamma)q^6+\cdots,
\end{aligned}
\end{equation}
where
\[
 \beta=\frac{|S(L)_{9/4}|}{2}\text{ and }
 \gamma=\frac{|S(L)_{17/4}|}{2}
\]
are nonnegative integers satisfying
\begin{equation}\label{eq:rank49-min4-basic-range}
\begin{aligned}
 &\beta\in\{0,1,\ldots,48\}\text{ and }
 &5392-71\beta\leq\gamma
 \leq\left\lfloor\frac{1473920+5419\beta}{3}\right\rfloor.
\end{aligned}
\end{equation}
Moreover, $\gamma\equiv\beta\pmod4$.
\end{prop}

\begin{proof}
Substituting the coefficients in Lemma~\ref{lem:theta-common-coefficients}
into~\eqref{eq:theta} gives~\eqref{eq:rank49-min4-theta}
and determines all its coefficients.
By Lemma~\ref{lem:theta-shadow-bounds} (iii), we have $\beta\leq49$,
and the Gram matrix of $S(L)_{9/4}\cap \mathcal C_1$ is $2I_\beta+J_\beta/4$.
If $\beta=49$, write $S(L)_{9/4}\cap \mathcal C_1=\{s_1,s_2,\ldots,s_{49}\}$.
The vectors $2s_1,2s_2,\ldots,2s_{49}$ generate
a sublattice $K$ of $L$ of full rank with Gram matrix $8I_{49}+J_{49}$.
Since this matrix has eigenvalues $8$ with multiplicity $48$ and $57$
with multiplicity $1$, we obtain
\[
 \det(K)=57\cdot8^{48}.
\]
However, since $L$ is unimodular, $\det(K)=|L:K|^2$ is a square in $\ZZ$,
while $57\cdot8^{48}$ is not a square.
This contradiction shows that $\beta\leq48$.

Since the coefficient of $q^6$ in $\theta_L(q)$ is nonnegative, we see that
\[
 48\gamma\leq23582720+86704\beta.
\]
This gives the upper bound on $\gamma$ in~\eqref{eq:rank49-min4-basic-range}.

It remains to prove the congruence and the lower bound on $\gamma$.
Let $N=|L_4|=71\beta+\gamma-5390$.

Expanding~\eqref{eq:theta} one term further, we find that the coefficient
of $q^7$ in $\theta_L(q)$ is
\[
 870912000-842400\beta-96\gamma.
\]
Since $48\gamma\leq23582720+86704\beta$ and $\beta\leq48$, we obtain
\begin{align*}
 870912000-842400\beta-96\gamma
 &\geq823746560-1015808\beta\\
 &\geq774987776>0.
\end{align*}
Hence, the coefficient of $q^7$ in $\theta_L(q)$ is positive,
and we may choose a vector $u\in L$ with $(u,u)=7$.
Applying~\eqref{eq:theta-common-weighted-relation} with $n=49$ to
$P_u(x)=(x,u)^2-\frac7{49}(x,x)$, we obtain
\begin{align*}
 &\sum_{s\in S(L)_{17/4}}(s,u)^2
  +47\sum_{s\in S(L)_{9/4}}(s,u)^2
  +4\sum_{y\in L_4}(y,u)^2\\
 &\qquad=\frac7{49}
  \left(\frac{17}{4}\cdot2\gamma
       +47\cdot\frac94\cdot2\beta+16N\right)\\
 &\qquad=\frac72(\gamma+55\beta-3520).
\end{align*}
Subtracting $(\gamma+47\beta)/2$, we obtain
\begin{align*}
 &\sum_{s\in S(L)_{17/4}}\left((s,u)^2-\frac14\right)
  +47\sum_{s\in S(L)_{9/4}}\left((s,u)^2-\frac14\right)
  +4\sum_{y\in L_4}(y,u)^2\\
 &\qquad=3\gamma+169\beta-12320.
\end{align*}
It follows from~\eqref{eq:theta-shadow-square-congruence} that
each of the first two sums is divisible by $4$.
The third term is also divisible by $4$, since $(y,u)\in\ZZ$ for $y\in L_4$.
Thus, we conclude that
\[
 3\gamma+169\beta-12320\equiv0\pmod4.
\]
This gives $\gamma\equiv\beta\pmod4$.
It follows that
\[
 |L_4|=71\beta+\gamma-5390\equiv72\beta-5390\equiv2\pmod4.
\]
Since $|L_4|$ is a nonnegative integer, we see that $|L_4|\geq2$,
which gives $\gamma\geq5392-71\beta$.
This completes the proof.
\end{proof}

\begin{cor}\label{cor:rank49-direct}
Every unimodular lattice in dimension $49$ has minimum norm at most $4$.
\end{cor}

\begin{proof}
Let $L$ be a unimodular lattice in dimension $49$.
Suppose that $\min(L)\geq4$.
By Proposition~\ref{prop:rank49-min4-theta},
the coefficient $71\beta+\gamma-5390$ of $q^4$ in $\theta_L(q)$
is at least $2$.
Hence, $L$ contains a vector of norm $4$.
Together with the assumption $\min(L)\geq4$, this gives $\min(L)=4$.
\end{proof}

Shi, Yang and Sol\'e~\cite{SYS} recently constructed a unimodular
lattice with minimum norm $4$ in dimension $49$.
Thus, we obtain the following:

\begin{cor}\label{cor:res49}
The largest minimum norm among unimodular lattices in dimension $49$ is $4$.
\end{cor}

\subsection{Dimension 50}

The upper bound $\mu(50)\leq5$ is given in~\cite[Table~2]{Ga04}.
We consider the theta series of unimodular lattices attaining this bound.

\begin{prop}\label{prop:rank50-min5-theta}
Let $L$ be a unimodular lattice in dimension $50$.
If $\min(L)=5$, then its theta series has one of the following forms:
\begin{align*}
\theta_{50,1}(q)&=1+252720q^5+20545200q^6+828921600q^7+\cdots,\\
\theta_{50,2}(q)&=1+250672q^5+20586160q^6+828593920q^7+\cdots\text{ and}\\
\theta_{50,3}(q)&=1+248624q^5+20627120q^6+828266240q^7+\cdots.
\end{align*}
Here $\beta=|S(L)_{5/2}|/2\in\{0,1,2\}$ and
$\theta_L(q)=\theta_{50,\beta+1}(q)$.
\end{prop}
\begin{proof}
Substituting the coefficients in Lemma~\ref{lem:theta-common-coefficients}
into~\eqref{eq:theta} and~\eqref{eq:theta-S}, we obtain
\begin{align*}
\theta_L(q)
&=1+(252720-2048\beta)q^5+(20545200+40960\beta)q^6\\
&\quad +(828921600-327680\beta)q^7+\cdots\text{ and}\\
\theta_{S(L)}(q)
&=2\beta q^{5/2}+(36400-140\beta)q^{9/2}\\
&\quad +(280130400+4730\beta)q^{13/2}+\cdots.
\end{align*}
By Lemma~\ref{lem:theta-shadow-bounds}~(i), we have
$\beta\in\{0,1,2\}$.
Substituting these values gives the three stated series.
\end{proof}

\subsection{Dimension 51}

The upper bound $\mu(51)\leq5$ is given in~\cite[Table~2]{Ga04}.
We consider the theta series of unimodular lattices attaining this bound.

\begin{prop}\label{prop:rank51-min5-theta}
Let $L$ be a  unimodular lattice in dimension $51$.
If $\min(L)=5$, then its theta series has the form
\begin{equation}\label{eq:rank51-min5-theta}
\begin{aligned}
\theta_L(q)
={}&1+(200124-1024\beta)q^5
+(17913240+18432\beta)q^6\\
 &+(781364880-122880\beta)q^7+\cdots,
\end{aligned}
\end{equation}
where
\[
 \beta=\frac{|S(L)_{11/4}|}{2}
\]
is a nonnegative integer satisfying
\begin{equation}\label{eq:rank51-min5-range}
 \beta\in\{0,1,\ldots,10,12\}.
\end{equation}
\end{prop}

\begin{proof}
Substituting the coefficients in Lemma~\ref{lem:theta-common-coefficients}
into~\eqref{eq:theta} gives~\eqref{eq:rank51-min5-theta}
and determines all its coefficients.
Also, \eqref{eq:theta-S} gives
\begin{align*}
 \theta_{S(L)}(q)
 &=2\beta q^{11/4}+(110160-138\beta)q^{19/4}\\
 &\quad +(641227760+4590\beta)q^{27/4}+\cdots.
\end{align*}
The restriction~\eqref{eq:rank51-min5-range} follows from
Lemma~\ref{lem:theta-shadow-bounds}~(ii).
\end{proof}

\subsection{Dimension 55}

Shi, Yang and Sol\'e~\cite{SYS} showed that $\mu(55)=5$.
We consider the theta series of optimal unimodular lattices
in dimension $55$.

\begin{prop}\label{prop:rank55-min5-theta}
Let $L$ be an optimal unimodular lattice in dimension $55$.
Then its theta series has the form
\begin{equation}\label{eq:rank55-min5-theta}
\begin{aligned}
 \theta_L(q)={}&1+(74492-5696\beta-64\gamma)q^5
  +(9373320+319104\beta+640\gamma)q^6\\
 &+(559096560-8912896\beta)q^7+\cdots,
\end{aligned}
\end{equation}
where
\[
 \beta=|S(L)_{7/4}|/2 \text{ and } \gamma=|S(L)_{15/4}|/2
\]
are nonnegative integers satisfying
\[
 \beta\in\{0,1\} \text{ and }
 0\leq\gamma\leq
 \begin{cases}
  1163 & \text{if } \beta=0,\\
  1042 & \text{if } \beta=1.
 \end{cases}
\]
\end{prop}

\begin{proof}
Since $\min(L)=5$, the coefficients of $q,q^2,q^3$ and $q^4$
in $\theta_L(q)$ are zero. Comparing coefficients in~\eqref{eq:theta} gives
\[
 (a_0,a_1,a_2,a_3,a_4)=(1,-110,3520,-32120,47410).
\]
Expanding~\eqref{eq:theta-S}, we obtain
\[
 \theta_{S(L)}(q)=\frac{a_6}{2^{17}}q^{7/4}
 +\left(-\frac{a_5}{32}-\frac{89a_6}{2^{17}}\right)q^{15/4}
 +\cdots.
\]
Hence $a_6=2^{18}\beta$ and $a_5=-64(89\beta+\gamma)$.
Substitution into~\eqref{eq:theta} gives~\eqref{eq:rank55-min5-theta}
and determines all its coefficients.

Suppose that $s,t\in S(L)_{7/4}$ and $t\neq\pm s$.
Since $S(L)$ is a coset of $L$ in $L_0^*$ and $2t\in L$,
both $s-t$ and $s+t$ are nonzero vectors of $L$. This gives
\[
 10\leq(s-t,s-t)+(s+t,s+t)=2(s,s)+2(t,t)=7,
\]
a contradiction. Thus, we conclude that $\beta\in\{0,1\}$.
If $\beta=0$, then $|L_5|=74492-64\gamma>0$, giving $\gamma\leq1163$.

Suppose now that $\beta=1$ and fix $s\in S(L)_{7/4}$.
For each vector $t\in S(L)_{15/4}$, the vectors $t-s$ and $t+s$
are nonzero vectors of $L$. Since $L$ is integral and $\min(L)=5$,
their norms are integers at least $5$. Their sum is
\[
 (t-s,t-s)+(t+s,t+s)=2(t,t)+2(s,s)=11.
\]
Hence, one of the two vectors has norm $5$ and the other has norm $6$.
This gives $(s,t)=\pm1/4$.
From each pair $\{t,-t\}$ choose the vector $t$ with $(s,t)=1/4$.
The $\gamma$ vectors $t-s$ have norm $5$ and are distinct.
Their negatives give a further $\gamma$ vectors, since
$(s,t-s)=-3/2$ for every chosen $t$. Therefore, we conclude that
\[
 2\gamma\leq |L_5|=68796-64\gamma,
\]
which gives $\gamma\leq1042$.
\end{proof}

\subsection{Dimension 57}

The upper bound $\mu(57)\leq5$ is given in~\cite[Table~2]{Ga04}.
We consider the theta series of unimodular lattices attaining this bound.

\begin{prop}\label{prop:rank57-min5-theta}
Let $L$ be a unimodular lattice in dimension $57$.
If $\min(L)=5$, then its theta series has the form
\begin{equation}\label{eq:rank57-min5-theta}
\begin{aligned}
 \theta_L(q)
 ={}&1+(49020-1392\beta-16\gamma)q^5
 +(6390080+73888\beta+96\gamma)q^6\\
 &+(446770560-1915968\beta+576\gamma)q^7+\cdots,
\end{aligned}
\end{equation}
where
\[
 \beta=\frac{|S(L)_{9/4}|}{2}\text{ and }
 \gamma=\frac{|S(L)_{17/4}|}{2}
\]
are nonnegative integers satisfying
\[
 \beta\in\{0,1\}\text{ and }0\leq\gamma\leq3063-87\beta.
\]
In particular, $|L_5|\equiv12\pmod{16}$.
\end{prop}

\begin{proof}
Substituting the coefficients in Lemma~\ref{lem:theta-common-coefficients}
into~\eqref{eq:theta} gives~\eqref{eq:rank57-min5-theta}
and determines all its coefficients.
To obtain a bound on $\beta$, suppose that $s,t\in S(L)_{9/4}$ and $t\neq\pm s$.
Since $S(L)$ is a coset of $L$ in $L_0^*$, we have $s-t\in L$.
Also, since $2t\in L$, we see that $s+t\in L$.
Since both $s-t$ and $s+t$ are nonzero, $\min(L)\geq5$ gives
\[
 10\leq(s-t,s-t)+(s+t,s+t)=2(s,s)+2(t,t)=9,
\]
a contradiction. Hence, $S(L)_{9/4}$ contains at most one pair
$\{s,-s\}$, which gives $\beta\in\{0,1\}$.

The coefficient of $q^5$ in~\eqref{eq:rank57-min5-theta} is
\[
 |L_5|=49020-16(87\beta+\gamma)\equiv12\pmod{16}.
\]
Since $|L_5|$ is nonnegative, it is at least $12$.
This gives $\gamma\leq3063-87\beta$.
\end{proof}

\subsection{Dimension 58}

The bound~\eqref{eq:RS} gives $\mu(58)\leq6$,
as also recorded in~\cite[Table~2]{Ga04}.

\begin{thm}\label{thm:common-upper-bounds}
Every unimodular lattice in dimension $58$ has minimum norm at most $5$.
\end{thm}

\begin{proof}
Let $L$ be a unimodular lattice in dimension $58$.
Suppose that $\min(L)\geq6$. 
If $s\in S(L)$ has norm $1/2$, then $2s\in L$ has norm $2$.
Hence, we conclude that $S(L)_{1/2}=\emptyset$.

The constant term and the vanishing of the coefficients of
$q,q^2,\ldots,q^5$ in~\eqref{eq:theta} give
\[
 (a_0,a_1,\ldots,a_5)=(1,-116,4060,-44080,96860,-44080).
\]
In~\eqref{eq:theta-S}, only the term with $j=7$ can contribute to the
coefficient of $q^{1/2}$, which is $-a_7/2^{26}$.
Since $S(L)_{1/2}=\emptyset$, we have $a_7=0$.
Then by expanding the terms with $j=6$ and $j=5$, we obtain
\[
 \theta_{S(L)}(q)
 =\frac{a_6}{2^{14}}q^{5/2}
  +\left(11020-\frac{86a_6}{2^{14}}\right)q^{9/2}+\cdots.
\]
Let $b=|S(L)_{5/2}|$. Comparing coefficients, we obtain
\begin{equation}\label{eq:common-shadow-counts}
 a_6=2^{14}b \text{ and } |S(L)_{9/2}|=11020-86b.
\end{equation}
Since $S(L)_{5/2}$ is a disjoint union of pairs $\{s,-s\}$, $b$ is even.
Substituting these coefficients into~\eqref{eq:theta}, we obtain
\[
 \theta_L(q)=1+(5380080+16384b)q^6
                 +(395884800-458752b)q^7+\cdots.
\]
From~\eqref{eq:common-shadow-counts}, we see that $11020-86b\geq0$.
This gives $b\leq128$.
Hence, the coefficient of $q^7$ is positive.
Therefore, we may choose
a vector $u\in L$ with $(u,u)=7$.

For the vector $u$ of norm $7$ chosen above, let
\[
 P_u(x)=(x,u)^2-\frac{7}{58}(x,x).
\]
This polynomial is homogeneous of degree $2$ and harmonic, since
\[
 \sum_{i=1}^{58}\frac{\partial^2 P_u}{\partial x_i^2}
 =2(u,u)-\frac{7}{58}\cdot116=14-14=0.
\]
Since $\min(L)\geq6$, we have $L_5=\emptyset$.
Applying~\eqref{eq:theta-common-weighted-relation} with $n=58$ to $P_u$ and
using~\eqref{eq:common-shadow-counts}, we obtain
\begin{align*}
 \sum_{s\in S(L)_{9/2}}(s,u)^2+62\sum_{s\in S(L)_{5/2}}(s,u)^2
 &=\frac{7}{58}
   \left(\frac92|S(L)_{9/2}|+62\cdot\frac52|S(L)_{5/2}|\right)\\
 &=\frac{7}{58}
   \left(\frac92(11020-86b)+62\cdot\frac52 b\right)\\
   &=7(855-4b).
\end{align*}
Also, \eqref{eq:common-shadow-counts} gives
\[
 \frac{|S(L)_{9/2}|+62|S(L)_{5/2}|}{4}
 =\frac{11020-86b+62b}{4}=2755-6b.
\]
Subtracting this quantity from the preceding sum, we obtain
\begin{equation}\label{eq:common-contradiction}
\begin{aligned}
 &\sum_{s\in S(L)_{9/2}}\left((s,u)^2-\frac14\right)
   +62\sum_{s\in S(L)_{5/2}}\left((s,u)^2-\frac14\right)\\
 &\qquad=7(855-4b)-(2755-6b)\\
 &\qquad\equiv2\pmod4,
\end{aligned}
\end{equation}
since $b$ is even.
On the other hand, $u$ has odd norm, so
\eqref{eq:theta-shadow-square-congruence} with $r=9/2$ and $r=5/2$
shows that each of the two sums on the left-hand side
of~\eqref{eq:common-contradiction} is divisible by $4$.
This is a contradiction. Therefore, we conclude that $\min(L)\leq5$.
\end{proof}

\begin{prop}\label{prop:rank58-min5-theta}
Let $L$ be a unimodular lattice in dimension $58$.
If $\min(L)=5$, then its theta series has the form
\begin{equation}\label{eq:rank58-min5-theta}
\begin{aligned}
\theta_L(q)
={}&1+(44080-688\beta-8\gamma)q^5
+(5203760+35520\beta+32\gamma)q^6\\
 &+(393768960-884480\beta+384\gamma)q^7+\cdots,
\end{aligned}
\end{equation}
where
\[
 \beta=\frac{|S(L)_{5/2}|}{2} \text{ and }
 \gamma=\frac{|S(L)_{9/2}|}{2}
\]
are nonnegative integers satisfying
\begin{equation}\label{eq:rank58-min5-basic-range}
 \beta\in\{0,1,2\} \text{ and } 0\leq\gamma\leq5509-86\beta.
\end{equation}
\end{prop}

\begin{proof}
Substituting the coefficients in Lemma~\ref{lem:theta-common-coefficients}
into~\eqref{eq:theta} gives~\eqref{eq:rank58-min5-theta}
and determines all its coefficients.
By Lemma~\ref{lem:theta-shadow-bounds}~(i), we have
$\beta\in\{0,1,2\}$.
Since the coefficient of $q^5$ in $\theta_L(q)$ is positive,
we have $\gamma<5510-86\beta$.
This gives~\eqref{eq:rank58-min5-basic-range}.
\end{proof}

Shi, Yang and Sol\'e~\cite{SYS} recently constructed a unimodular
lattice with minimum norm $5$ in dimension $58$.
Thus, we obtain the following:

\begin{cor}\label{cor:res58}
The largest minimum norm among unimodular lattices in dimension $58$ is $5$.
\end{cor}

%

\subsection{Dimension 59}

Shi, Yang and Sol\'e~\cite{SYS} showed that $\mu(59)=5$.
We consider the theta series of optimal unimodular lattices
in dimension $59$.

\begin{prop}\label{prop:rank59-min5-theta}
Let $L$ be an optimal unimodular lattice in dimension $59$.
Then its theta series has the form
\begin{equation}\label{eq:rank59-min5-theta}
\begin{aligned}
 \theta_L(q)={}&1+(44604-340\beta-4\gamma)q^5
  +(4212600+17064\beta+8\gamma)q^6\\
 &+(343705680-406944\beta+224\gamma)q^7+\cdots,
\end{aligned}
\end{equation}
where 
\[
\beta=|S(L)_{11/4}|/2 \text{ and }\gamma=|S(L)_{19/4}|/2
\]
are nonnegative integers satisfying
\[
 \beta\in\{0,1,\ldots,10,12\},
 \beta+\gamma\equiv0\pmod4 \text{ and }
 0\leq\gamma\leq11148-85\beta.
\]
\end{prop}

\begin{proof}
Substituting the coefficients in Lemma~\ref{lem:theta-common-coefficients}
into~\eqref{eq:theta} gives~\eqref{eq:rank59-min5-theta}
and determines all its coefficients.
Lemma~\ref{lem:theta-shadow-bounds}~(ii) gives
the stated possibilities for $\beta$.

Choose a vector $w\in S(L)$ and let $u=2w$ and $t=(u,u)$.
By the description of the shadow in~\cite[p.~1090]{Pache05},
$u$ is a characteristic vector of $L$, that is,
\[
  (u,x)\equiv(x,x)\pmod 2
\]
for every vector $x\in L$.
By Proposition~3 (i) in~\cite{Pache05}, we have
\[
  t=4(w,w)\equiv59\equiv3\pmod 8.
\]

The coefficient of $q^5$ in~\eqref{eq:rank59-min5-theta} gives
\[
  |L_5|=44604-340\beta-4\gamma.
\]
As in the proof of Theorem~\ref{thm:common-upper-bounds},
the following polynomial
\[
  P(x)=(x,u)^2-\frac{t}{59}(x,x)
\]
is homogeneous of degree $2$ and harmonic.
Applying~\eqref{eq:theta-common-weighted-relation} with $n=59$
to $P$ and multiplying both sides by $2$, we obtain
\begin{align*}
 &2\sum_{s\in S(L)_{19/4}}(s,u)^2
  +122\sum_{s\in S(L)_{11/4}}(s,u)^2
  -2\sum_{y\in L_5}(y,u)^2\\
 &\qquad=\frac{2t}{59}
   \left(\frac{19}{2}\gamma
         +61\cdot\frac{11}{2}\beta-5|L_5|\right)\\
 &\qquad=t(69\beta+\gamma-7560).
\end{align*}

Since $u$ has odd norm, \eqref{eq:theta-shadow-square-congruence} gives
\[
  2\sum_{s\in S(L)_r}(s,u)^2
  \equiv\frac{|S(L)_r|}{2}\pmod 8
\]
for $r=19/4$ and $11/4$.
Thus, the two sums over the shadow, with their respective factors
$2$ and $122$, are congruent to $\gamma$ and $61\beta$ modulo $8$.
Since $u$ is characteristic, $(y,u)$ is odd for every vector $y\in L_5$.
Hence, we see that $(y,u)^2\equiv1\pmod 8$.
The above expression for $|L_5|$ gives $|L_5|\equiv0\pmod 4$, so
\[
  -2\sum_{y\in L_5}(y,u)^2\equiv-2|L_5|\equiv0\pmod 8.
\]
Using $t\equiv3\pmod 8$, we obtain
\[
  \gamma+61\beta\equiv3(69\beta+\gamma-7560)\pmod 8.
\]
This gives $2(\beta+\gamma)\equiv0\pmod 8$, and hence
$\beta+\gamma\equiv0\pmod 4$.

Finally, since $\min(L)=5$, we have
\[
  |L_5|=4(11151-85\beta-\gamma)>0.
\]
Thus, $85\beta+\gamma<11151$.
Since $85\beta+\gamma\equiv\beta+\gamma\equiv0\pmod 4$,
we obtain $85\beta+\gamma\leq11148$.
This gives $\gamma\leq11148-85\beta$.
\end{proof}

\subsection{Dimensions 52, 53, 54 and 56}

It is known that $\mu(52)=\mu(53)=\mu(54)=5$ and $\mu(56)=6$
(see~\cite{Ga04}, \cite{NG} and~\cite{SYS}).
In dimension $52$, the theta series of an optimal unimodular lattice and
of its shadow were given in~\cite[Section~3]{Ha20}. It was also shown
there that such a lattice is $s$-extremal if and only if its kissing
number is $157248$.
In dimension $56$, even unimodular lattices with minimum norm $6$
are known~\cite{NG}. The coefficient of $q^6$ in the theta series of
every such lattice is $15590400$ (see~\cite[Chapter~7, Section~7]{SPLAG}),
and all odd powers of $q$ have coefficient zero.
We give necessary conditions on the theta series of an odd unimodular
lattice attaining $\mu(n)$ in each of these four dimensions.

\begin{prop}\label{prop:remaining-theta}
Let $L$ be an optimal odd unimodular lattice in dimension $n$, where
$n\in\{52,53,54,56\}$. Then its theta series has the form
\[
 \theta_L(q)=1+A_5q^5+A_6q^6+A_7q^7+\cdots,
\]
where $A_5$, $A_6$ and $A_7$ are given in Table~\ref{Tab:remaining-theta}.
The parameters in the table are nonnegative integers satisfying
the following conditions.
\begin{enumerate}
\item If $n=52$, then $\beta=|S(L)_3|/2$ and $\beta\leq72$.
\item If $n=53$ or $54$, let $r=(n-48)/4$ and put
$\beta=|S(L)_r|/2$ and $\gamma=|S(L)_{r+2}|/2$.
Then $\beta\in\{0,1\}$ and
\[
  0\leq\gamma\leq
  \begin{cases}
    479(1-\beta) & \text{if } n=53,\\
    745-109\beta & \text{if } n=54.
  \end{cases}
\]
\item If $n=56$, then $\beta=|S(L)_2|/2\in\{0,1\}$.
\end{enumerate}
\end{prop}

\begin{table}[thb]
\caption{Coefficients of $\theta_L(q)$ for optimal odd unimodular lattices}
\label{Tab:remaining-theta}
\centering
\medskip
{\small
\begin{tabular}{c|l|l|l}
\noalign{\hrule height1pt}
$n$ & \multicolumn{1}{c|}{$A_5$} & \multicolumn{1}{c|}{$A_6$} 
& \multicolumn{1}{c}{$A_7$}\\
\hline
$52$ & $157248-512\beta$ & $15462720+8192\beta$ & $729181440-43008\beta$\\
$53$ & $122748-23296\beta$ & $13216080+1374720\beta$ & $673829280-41057280\beta$\\
     & $-256\gamma$ & $+3584\gamma$ & $-13312\gamma$\\
$54$ & $95472-11520\beta$ & $11185200+662528\beta$ & $616705920-19150848\beta$\\
     & $-128\gamma$ & $+1536\gamma$ & $-3072\gamma$\\
$56$ & $0$ & $8250368+131072\beta$ & $503316480-4194304\beta$\\
\noalign{\hrule height1pt}
\end{tabular}
}
\end{table}

\begin{proof}
The constant term in~\eqref{eq:theta} gives $a_0=1$.
Since the coefficient of $q^i$ is zero for
$1\leq i<\min(L)$, comparing coefficients successively determines
$a_1,a_2,\ldots,a_{\min(L)-1}$. In~\eqref{eq:theta-S}, the term indexed
by $j$ has the expansion
\[
 (-1)^j2^{n-12j}a_j q^{(n-8j)/4}
 (1+(n-24j)q^2+\cdots).
\]
For $n=52$, since a vector $s\in S(L)_1$ gives a vector $2s\in L$
of norm $4$, $S(L)_1$ is empty. 
For $n=56$, the constant term of $\theta_{S(L)}(q)$ is zero,
since the zero vector is the only vector $\0$ of norm $0$ and belongs to $L$.
In the preceding expansion, the exponent $(56-8j)/4=14-2j$
is positive for $0\leq j\leq6$ and is zero for $j=7$.
Thus, only the term with $j=7$ can contribute to the constant term.
Comparing constant terms gives $0=-a_7/16^7$.  This gives  $a_7=0$.
Comparing coefficients in the shadow expansion and using the values
already determined from~\eqref{eq:theta}, we obtain
\[
 (a_5,a_6)=
 \begin{cases}
 (-512\beta,0) & \text{if }n=52,\\
 (-256(91\beta+\gamma),2^{20}\beta) & \text{if }n=53,\\
 (-128(90\beta+\gamma),2^{19}\beta) & \text{if }n=54,\\
 (-59136,2^{17}\beta) & \text{if }n=56.
 \end{cases}
\]
Substitution into~\eqref{eq:theta} gives Table~\ref{Tab:remaining-theta}
and determines all the coefficients of $\theta_L(q)$.

Since $S(L)$ is a coset of $L$ in $L_0^*$ and $2t\in L$ for
$t\in S(L)$, we have
$s-t,s+t\in L$ for any vectors $s,t\in S(L)$. If $t\ne\pm s$, then we obtain
\begin{equation}\label{eq:remaining-shadow-sum}
 2\min(L)\leq (s-t,s-t)+(s+t,s+t)=2(s,s)+2(t,t).
\end{equation}
In particular, for $0<r<\min(L)/2$, there is at most one pair
$\{s,-s\}$ in $S(L)_r$. This gives $\beta\in\{0,1\}$ for
$n\in\{53,54,56\}$.

Suppose that $n=52$. For $t\in S(L)_3$, the vector $2t$ has even norm
$12$ and belongs to $L_0$. Hence, $t$ and $-t$ belong to the same shadow
coset. If $s,t\in S(L)_3$ lie in the same coset and $t\ne\pm s$,
then $s-t,s+t\in L_0$. Their norms are even, at least $6$, and sum to $12$.
Thus, we see that $(s,t)=0$. If $s$ and $t$ lie in different shadow cosets, then
$s-t$ and $s+t$ lie in $L\setminus L_0$. Their norms are odd, at least
$5$, and sum to $12$. Hence, their norms are $5$ and $7$, giving
$(s,t)=\pm1/2$. If only one coset contains vectors of norm $3$,
orthogonality gives $\beta\leq52$. If both do, choose one vector
$s_i$ from each of the $k$ pairs in one coset and fix $t$ of norm $3$
in the other. 
Since $(s_i,s_i)=3$ for each $i$ and $(s_i,s_j)=0$ for $i\ne j$,
the vectors $s_1/\sqrt{3},s_2/\sqrt{3},\ldots,s_k/\sqrt{3}$ are orthonormal.
Applying Bessel's inequality to $t$ and these vectors, we obtain
\[
  \sum_{i=1}^k \left(t,\frac{s_i}{\sqrt{3}}\right)^2
  \leq (t,t)=3.
\]
Since $(t,s_i)=\pm1/2$ for each $i$, every term in the sum is $1/12$.
Thus $k/12\leq3$, which gives $k\leq36$.
Interchanging the roles of the two shadow cosets, we see that the
other coset also contains at most $36$ pairs of vectors of norm $3$.
Hence, we conclude that $\beta\leq36+36=72$.

For $n=53$ and $54$, the bounds for $\beta=0$ follow from
$A_5=|L_5|>0$ in Table~\ref{Tab:remaining-theta}.
Suppose that $\beta=1$ and fix a vector $s\in S(L)_r$.
For $n=53$, a vector $t\in S(L)_{13/4}$ gives
$10\leq2(5/4)+2(13/4)=9$ by~\eqref{eq:remaining-shadow-sum}.
Thus, we have $\gamma=0$. For $n=54$ and $t\in S(L)_{7/2}$, 
since the norms of
$t-s$ and $t+s$ sum to $10$,  both norms are $5$ and $(s,t)=0$.
As $t$ varies, the $2\gamma$ vectors $t+s$ are distinct, as are the
$2\gamma$ vectors $t-s$. These two sets are disjoint, since their inner
products with $s$ are $3/2$ and $-3/2$, respectively. Therefore, we conclude that
\[
 4\gamma\leq|L_5|=83952-128\gamma,
\]
which gives $\gamma\leq636$.
\end{proof}

We give the parameters of some known optimal unimodular lattices
in the notation of Proposition~\ref{prop:remaining-theta}.

\begin{itemize}
\item Dimension $52$:
The shadow of the lattice constructed in~\cite{Ga04} contains
$104$ vectors of norm $3$ (see~\cite[Section~3]{Ha20}).
Hence, we obtain $\beta=104/2=52$.
The three nonisomorphic $s$-extremal lattices constructed
in~\cite{Ha20} have $\beta=0$.

\item Dimension $53$:
Let $L$ be the lattice constructed in~\cite[Table~1]{SYS}.
Using \textsc{Magma}~\cite{Magma}, we counted the vectors of norm $5$
in $L$ and obtained $|L_5|=122748$.
By Proposition~\ref{prop:remaining-theta}, we have
\[
  122748=|L_5|=122748-23296\beta-256\gamma.
\]
Since $\beta$ and $\gamma$ are nonnegative, both are zero.
Thus, we see that $a_5=a_6=0$, and comparing coefficients through $q^4$
in~\eqref{eq:theta} gives $a_4=27030$.
Therefore, the leading term in~\eqref{eq:theta-S} comes from $j=4$
and is $2^5\cdot27030q^{21/4}=864960q^{21/4}$.
Hence, we have $\min(S(L))=21/4$ and
\[
  8\min(L)+4\min(S(L))=8\cdot5+4\cdot\frac{21}{4}=53+8.
\]
Thus, $L$ is $s$-extremal~\cite[Theorem~2.1 and Definition~2.2]{Ga07}.

\item Dimension $54$:
The lattice constructed in~\cite[Section~4.2.1]{Ga04} has kissing
number $95472$ and its shadow has minimum norm $11/2$.
Hence, we conclude that $(\beta,\gamma)=(0,0)$ in our notation.
This lattice is $s$-extremal~\cite[Section~2.3]{Ga07}.
\end{itemize}

\section{Construction of optimal unimodular lattices}\label{Sec:C}

In this section, we use Construction~A to construct unimodular lattices
with minimum norm $4$ in dimension $49$ and
unimodular lattices with minimum norm $5$ in dimensions 
$50,51,55,57,58$ and $59$.
Together with the upper bounds in Section~\ref{Sec:T}, our constructions
give $\mu(50)=\mu(51)=\mu(57)=5$.
We used \textsc{Magma}~\cite{Magma} to search for self-dual codes
and to verify the minimum norms of the resulting lattices.
To determine the theta series and kissing numbers, we used
the formulas in Section~\ref{Sec:T} together with counts of vectors
of prescribed norms in the lattices or their shadows.
Throughout this section, 
let $\tau$ denote the kissing number of the lattice under consideration.

\subsection{Matrices used in the constructions}

An $n \times n$ circulant matrix has the following form:
\[
\left(
\begin{array}{cccccc}
r_0&r_1 &\cdots      &r_{n-2} &r_{n-1} \\
r_{n-1}&r_0& \cdots &r_{n-3} &r_{n-2} \\
\vdots & \ddots & \ddots&& \vdots\\
r_2 &        &\ddots & \ddots& r_1\\
r_1&r_2& \cdots   &r_{n-1}&r_0 \\
\end{array}
\right).
\]
Consider the following matrix:
\begin{equation} \label{eq:GM}
\left(
\begin{array}{cc}
I_n &  R 
\end{array}
\right),
\end{equation}
where $I_n$ denotes the identity matrix of order $n$
and $R$ is an $n \times n$ circulant matrix.
The codes generated by the rows of matrices of the form~\eqref{eq:GM}
are called \emph{double circulant codes}.

For the constructions in odd dimensions, we also use matrices
of the following form:
\begin{equation}\label{eq:odd-matrix}
\left(
\begin{array}{ccccccccc}
 & &  & a &  & &  \\ 
  & I_{n} & & \vdots && R &\\
 & &  & a &  & &  \\ 
0&\cdots & 0 & b & c & \cdots & c
\end{array}
\right),
\end{equation}
where $R$ is an $n\times n$ circulant matrix.
The codes generated by the rows of these matrices have length $2n+1$.
We call such codes \emph{modified double circulant codes}.
The form in~\eqref{eq:odd-matrix} was inspired by the matrix used
in~\cite[p.~159]{GH} to construct a self-dual $\ZZ_4$-code $C_{39}$
of length $39$. The corresponding unimodular lattice $A_4(C_{39})$ is optimal.

\subsection{Dimension 50}

We found a double circulant
self-dual $\ZZ_5$-code $C_{50}$ of length $50$ such that
$A_5(C_{50})$ has minimum norm $5$.
The code $C_{50}$ is generated by the rows of a matrix of the form~\eqref{eq:GM},
where the first row of the circulant matrix $R$ is
\[
(2,3,4,3,2,0,1,4,1,0,1,0,4,0,0,2,2,4,1,3,2,2,0,4,2).
\]

Combining this construction with the upper bound $\mu(50)\leq5$
given in~\cite[Table~2]{Ga04}, we obtain the following:

\begin{thm}\label{thm:L50}
\begin{enumerate}
\item
The largest minimum norm among unimodular lattices in dimension $50$
is $5$.
\item
The lattice $A_5(C_{50})$ is optimal and has minimum norm $5$.
\end{enumerate}
\end{thm}

We determine the theta series of every optimal unimodular lattice
$A_5(C)$ arising from a double circulant self-dual $\ZZ_5$-code $C$
of length $50$.

\begin{prop}\label{cor:rank50-double-circulant}
Let $C$ be a double circulant self-dual $\ZZ_5$-code of length $50$.
If $\min(A_5(C))=5$, then
$A_5(C)$ has theta series $\theta_{50,1}(q)$.
\end{prop}

\begin{proof}
By Proposition~\ref{prop:rank50-min5-theta}, the possible kissing
numbers of $A_5(C)$ are $252720$, $250672$ and $248624$.
Of these, only $252720$ is divisible by $5$.
Therefore, it suffices to show that the kissing number of $A_5(C)$
is divisible by $5$.

Suppose that $C$ is generated by the rows of the matrix
$G=\left(
\begin{array}{cc}
I_{25} &  R 
\end{array}
\right)$,
where $R$ is a circulant matrix.
Let $\sigma$ be the coordinate permutation:
\[
  \sigma=(1,2,\ldots,25)(26,27,\ldots,50).
\]
Since $R$ is circulant, $\sigma$ permutes
the rows of $G$.
Thus, $\sigma$ preserves $C$.
Hence, $\sigma$ induces an automorphism of $A_5(C)$ of order $25$.

Suppose that $\sigma$ fixes a vector $x=(x_1,x_2,\ldots,x_{50})\in A_5(C)$
of norm $5$.
It follows from the form of $\sigma$ that
\[
 x_1=x_2=\cdots=x_{25} \text{ and }
 x_{26}=x_{27}=\cdots=x_{50}.
\]
Since $\sqrt{5}\,x\in\ZZ^{50}$, we may write
\[
 x=\frac{1}{\sqrt{5}}
 \bigl(\underbrace{a,a,\ldots,a}_{25},
       \underbrace{b,b,\ldots,b}_{25}\bigr)
\]
for $a,b\in\ZZ$.
Then $(x,x)=5(a^2+b^2)=5$, which gives $a^2+b^2=1$.
This gives
\[
  (a,b)=(\pm1,0) \text{ or } (0,\pm1).
\]
By the definition of $A_5(C)$,
the reduction of $(a,a,\ldots,a,b,b,\ldots,b)$ modulo $5$ is a codeword of $C$
with exactly one nonzero block.
However, every codeword of $C$ has the form $(u,uR)$, and
the self-duality of $C$ gives $RR^T=4I_{25}$,
where $R^T$ denotes the transpose of $R$.
Hence, $R$ is nonsingular.
This gives $uR=\0$ if and only if $u=\0$.
This is a contradiction.
Therefore, every orbit of $\langle\sigma\rangle$ on the vectors of norm
$5$ in $A_5(C)$ has size $5$ or $25$. It follows that the kissing number
of $A_5(C)$ is divisible by $5$.
By Proposition~\ref{prop:rank50-min5-theta}, we conclude that
$A_5(C)$ has theta series $\theta_{50,1}(q)$.
\end{proof}

If $C$ is a double circulant self-dual $\ZZ_5$-code of length $50$
and $A_5(C)$ has minimum norm $5$, then $\beta=0$ by the preceding
proposition, and the shadow expansion gives $\min(S(A_5(C)))=9/2$.
Thus, we have $8\min(A_5(C))+4\min(S(A_5(C)))=50+8$, and $A_5(C)$ is
$s$-extremal~\cite[Theorem~2.1 and Definition~2.2]{Ga07}.

\begin{prop}\label{prop:rank50-s-extremal}
There exists an $s$-extremal unimodular lattice with minimum norm $5$ in dimension $50$.
\end{prop}

\subsection{Dimensions 51 and 49}

We found a modified double circulant self-dual $\ZZ_9$-code
$C_{51}$ of length $51$ such that $A_9(C_{51})$ has minimum norm $5$.
The code $C_{51}$ is generated by the rows of a matrix of the
form~\eqref{eq:odd-matrix}, where $(a,b,c)=(1,3,3)$ and
the first row of the $25 \times 25$ circulant matrix $R$ is
\[
(4,2,1,7,3,7,3,4,3,8,5,6,2,2,3,1,7,0,7,6,6,3,5,7,5).
\]

Combining this construction with the upper bound $\mu(51)\leq5$
given in~\cite[Table~2]{Ga04}, we obtain the following:

\begin{thm}\label{thm:L51}
\begin{enumerate}
\item
The largest minimum norm  among unimodular lattices in dimension $51$
is $5$.
\item
The lattice $A_9(C_{51})$ is optimal and has minimum norm $5$.
\end{enumerate}
\end{thm}

We verified that the shadow of
$A_9(C_{51})$ contains exactly two vectors of norm $11/4$, namely
\[
  \pm(\underbrace{1,\ldots,1}_{25},7,\underbrace{-1,\ldots,-1}_{25})/6.
\]
Thus, we see that $\beta=1$
in Proposition~\ref{prop:rank51-min5-theta} and
$A_9(C_{51})$ has the following theta series:
\[
1+199100q^5+17931672q^6+781242000q^7+\cdots.
\]

Now we construct nonisomorphic 
optimal unimodular lattices in dimension $49$ from $A_9(C_{51})$.
Let $D=\{0,3,6\} \subset \ZZ_9$ and let $1\leq i<j\leq51$.
Let $\pi_{i,j}$ denote the map that deletes coordinates $i$ and $j$.
By subtracting the self-dual code $D^2$ from $C_{51}$ on these
two coordinates, we obtain the following code:
\[
  C_{49}(i,j)=
  \{\pi_{i,j}(c)\mid c=(c_1,c_2,\ldots,c_{51})\in C_{51},
                          c_i,c_j\in D\}.
\]
To see that $C_{49}(i,j)$ is self-dual, move coordinates $i$ and $j$
to the last two positions, keeping the other coordinates in their
original order, and denote the resulting code by $C$.
Let $U=\{(x,y) \mid x \in \ZZ_9^{49}, y\in D^2\}$. Since $D^2$ is self-dual,
we see that $U^\perp=\{(\0,y) \mid y \in D^2\}$. 
Then we obtain
\[
  (C\cap U)^\perp=C^\perp+U^\perp=C+U^\perp.
\]
Thus, for every vector $x\in \ZZ_9^{49}$, we obtain
\begin{align*}
  x\in C_{49}(i,j)^\perp
  &\iff (x,\0)\in(C\cap U)^\perp\\
  &\iff (x,d)\in C \text{ for some } d\in D^2\\
  &\iff x\in C_{49}(i,j).
\end{align*}
We found $24$
unimodular lattices $A_9(C_{49}(i,j))$ with minimum norm $4$ in dimension $49$,
where the pairs  $(i,j)$ are listed in Table~\ref{Tab:49}.
By Proposition~\ref{prop:rank49-min4-theta},
the theta series of each $A_9(C_{49}(i,j))$ is determined by
the parameters $\beta$ and $\gamma$ defined there.
For each pair $(i,j)$, Table~\ref{Tab:49} lists the triple
$(\beta,\gamma,\tau)$.
The $24$ lattices have distinct theta series.
The lattice $A_9(C_{49}(1,2))$, which has kissing number $1418$,
is listed as \texttt{Uni49min4} in the catalogue~\cite{NG}.
For the lattice in dimension $49$ constructed in~\cite[Table~1]{SYS},
we obtained $(\beta,\gamma,\tau)=(0,7448,2058)$.
Since these $25$ lattices have distinct kissing numbers, 
they are nonisomorphic.

\begin{prop}\label{prop:L49}
There exist at least $25$ nonisomorphic optimal unimodular lattices
with minimum norm $4$ in dimension $49$.
\end{prop}

\begin{table}[thb]
\caption{Unimodular lattices $A_9(C_{49}(i,j))$}
\label{Tab:49}
\centering
\medskip
{\small
\begin{tabular}{c|c||c|c}
\noalign{\hrule height1pt}
$(i,j)$& $(\beta,\gamma,\tau)$&
$(i,j)$& $(\beta,\gamma,\tau)$\\
\hline
$(1,2)$  & $(0,6808,1418)$ & $(1,29)$  & $(0,6776,1386)$ \\
$(1,3)$  & $(0,6816,1426)$ & $(1,31)$  & $(0,6904,1514)$ \\
$(1,4)$  & $(0,6792,1402)$ & $(1,32)$  & $(0,6688,1298)$ \\
$(1,5)$  & $(0,6704,1314)$ & $(1,33)$  & $(0,6768,1378)$ \\
$(1,6)$  & $(0,6800,1410)$ & $(1,34)$  & $(0,6888,1498)$ \\
$(1,7)$  & $(0,6728,1338)$ & $(1,35)$  & $(0,6936,1546)$ \\
$(1,9)$  & $(0,6744,1354)$ & $(1,38)$  & $(0,6856,1466)$ \\
$(1,10)$ & $(0,6784,1394)$ & $(1,41)$  & $(0,6736,1346)$ \\
$(1,11)$ & $(0,6848,1458)$ & $(1,48)$  & $(0,6872,1482)$ \\
$(1,13)$ & $(0,6824,1434)$ & $(27,30)$ & $(0,6760,1370)$ \\
$(1,26)$ & $(0,6448,1058)$ & $(27,31)$ & $(0,6720,1330)$ \\
$(1,27)$ & $(0,6832,1442)$ & $(27,37)$ & $(0,6672,1282)$ \\
\noalign{\hrule height1pt}
\end{tabular}
}
\end{table}

\subsection{Dimension 55}

We found a modified double circulant self-dual $\ZZ_{25}$-code
$C_{55}$ of length $55$ such that $A_{25}(C_{55})$ has minimum norm $5$.
The code $C_{55}$ is generated by the rows of a matrix of the
form~\eqref{eq:odd-matrix}, where $(a,b,c)=(2,5,10)$ and
the first row of the $27 \times 27$ circulant matrix $R$ is
\[
 \begin{aligned}
 &(5,10,3,7,17,6,9,13,22,10,18,19,6,11,\\
 &\qquad 5,21,14,14,20,23,12,11,11,21,17,22,7).
 \end{aligned}
\]

The shadow of $A_{25}(C_{55})$ contains exactly two vectors of norm $7/4$,
namely
\[
 \pm(\underbrace{1,\ldots,1}_{27},-11,\underbrace{-1,\ldots,-1}_{27})/10.
\]
We also verified that the shadow contains no vectors of norm $15/4$.
Thus, by Proposition~\ref{prop:rank55-min5-theta},
we obtained $(\beta,\gamma,\tau)=(1,0,68796)$.
For the lattice in dimension $55$ constructed in~\cite[Table~1]{SYS},
we obtained $(\beta,\gamma,\tau)=(0,0,74492)$.
Since these two lattices have distinct kissing numbers, they are nonisomorphic.

\begin{prop}\label{prop:L55}
There exist at least two nonisomorphic optimal unimodular lattices
with minimum norm $5$ in dimension $55$.
\end{prop}

\subsection{Dimension 57}

We found $3$ modified double circulant self-dual $\ZZ_{25}$-codes
$C_{57,i}$ of length $57$ such that $A_{25}(C_{57,i})$ has
minimum norm $5$ for $i=1,2,3$.
The codes $C_{57,i}$ are generated by the rows of matrices of the
form~\eqref{eq:odd-matrix}, where $(a,b,c)=(1,5,20)$ and
the first rows of the $28 \times 28$ circulant matrices $R$ are
\begin{align*}
&(6,6,17,22,9,6,15,2,4,20,12,7,17,4,\\
      &\qquad 21,15,11,8,13,9,16,17,16,9,19,24,8,3),\\
&(16,19,10,7,17,14,9,22,21,24,19,2,15,4,\\
      &\qquad 15,4,3,0,2,18,7,15,16,3,8,9,1,11) \text{ and}\\
&(1,5,24,24,21,20,7,0,16,3,12,10,13,11,\\
      &\qquad 23,10,20,5,23,6,20,4,20,15,21,19,13,20),
\end{align*}
respectively.
By Proposition~\ref{prop:rank57-min5-theta},
the theta series of each $A_{25}(C_{57,i})$ is determined by
the parameters $\beta$ and $\gamma$ defined there.
For \(i=1,2,3\), the triples \((\beta,\gamma,\tau)\) are
\begin{align*}
(1,322,42476), (1,350,42028)\text{ and }(1,490,39788),
\end{align*}
respectively.
Since these $3$ lattices have distinct kissing numbers, they are
nonisomorphic.
Combining these constructions with the upper bound $\mu(57)\leq5$
given in~\cite[Table~2]{Ga04}, we obtain the following:

\begin{thm}\label{thm:L57}
\begin{enumerate}
\item
The largest minimum norm among unimodular lattices in dimension $57$
is $5$.
\item
The $3$ lattices $A_{25}(C_{57,i})$ $(i=1,2,3)$ are 
nonisomorphic optimal unimodular lattices with minimum norm $5$.
\end{enumerate}
\end{thm}

\subsection{Dimension 58}

We found $6$ double circulant self-dual $\ZZ_5$-codes $C_{58,i}$ of length $58$
such that $A_5(C_{58,i})$ has minimum norm $5$ for $i=1,2,\ldots,6$.
Each code $C_{58,i}$ is generated by the rows of a matrix
of the form~\eqref{eq:GM}.

\begin{table}[thb]
\caption{Unimodular lattices $A_5(C_{58,i})$}
\label{Tab:58}
\centering
\medskip
{\footnotesize
\begin{tabular}{c|c|c}
\noalign{\hrule height1pt}
$i$ & $r$ & $(\beta,\gamma,\tau)$  \\
\hline
1&$(3,3,2,2,4,1,0,4,1,1,2,2,1,4,3,2,2,4,4,3,2,2,3,2,1,1,3,3,2)$&$(0,2030,27840)$\\
2&$(1,4,4,4,2,4,3,2,4,4,4,3,4,3,2,0,4,4,3,4,1,2,1,1,2,4,3,3,2)$&$(0,1798,29696)$\\
3&$(1,2,4,0,2,2,0,1,4,1,0,0,0,2,2,3,2,2,1,0,0,1,2,1,4,3,0,3,4)$&$(0,1682,30624)$\\
4&$(2,2,3,3,4,3,4,0,2,1,1,4,2,0,0,0,1,1,0,2,0,1,4,1,4,4,0,1,2)$&$(0,1508,32016)$\\
5&$(0,4,3,3,4,3,1,0,3,0,1,4,1,0,3,3,1,0,2,3,0,4,1,0,3,2,2,1,0)$&$(0,1450,32480)$\\
6&$(0,1,1,1,1,0,3,1,4,1,0,2,1,0,1,0,3,1,1,4,0,3,4,3,4,4,2,4,2)$&$(0,1334,33408)$\\
\noalign{\hrule height1pt}
\end{tabular}
}
\end{table}

By Proposition~\ref{prop:rank58-min5-theta},
the theta series of each $A_5(C_{58,i})$ is determined by
the parameters $\beta$ and $\gamma$ defined there.
For each $i$, Table~\ref{Tab:58} lists the first row $r$ of $R$
and the triple $(\beta,\gamma,\tau)$.
For the lattice in dimension $58$ constructed in~\cite[Table~2]{SYS},
we obtained $(\beta,\gamma,\tau)=(0,1160,34800)$.
Since these $7$ lattices have distinct kissing numbers, 
they are nonisomorphic.

\begin{prop}\label{prop:L58}
There exist at least $7$ nonisomorphic  optimal unimodular lattices with minimum norm $5$ in dimension $58$.
\end{prop}

\subsection{Dimension 59}

We describe the construction of the lattice listed as
\texttt{Uni59min5} in the catalogue~\cite{NG}.
Let $\chi$ denote the quadratic character of $\FF_{59}$ with $\chi(0)=0$,
where $\FF_{59}$ denotes the finite field of order $59$.
Let $Q=(\chi(j-i))_{i,j\in\FF_{59}}$ and let
\[
  H=
\left(
\begin{array}{cccc}
0 & 1 & \cdots & 1 \\
-1 &  & & \\
\vdots &  & Q& \\
-1 &  & & \\
\end{array}
\right).
\]
The rows of $2I_{60}+H$ generate a self-dual $\ZZ_9$-code $C_{60}$.
Using the coordinates given in~\cite{NG}, we verified that $A_9(C_{60})$
is isomorphic to the extremal odd unimodular lattice $P_{60q}$
constructed in~\cite{Ga04}.
Let $\pi_1$ denote the map that deletes the first coordinate.
By subtracting the self-dual code $D=\{0,3,6\} \subset \ZZ_9$ from $C_{60}$
on the first coordinate, we obtain the $\ZZ_9$-code
\[
  C_{59,1}=\{\pi_1(c)\mid c=(c_1,c_2,\ldots,c_{60})\in C_{60},
                         c_1\in D\}.
\]
The code $C_{59,1}$ is self-dual and is generated by the rows
of the circulant matrix $2I_{59}+Q-J_{59}$.
Every nonzero vector of $A_9(C_{59,1})$ can be lifted to a vector
of $A_9(C_{60})$ whose first coordinate is $0$, $1$ or $-1$.
Since $\min(A_9(C_{60}))=6$, deleting this coordinate gives
$\min(A_9(C_{59,1}))\geq5$. Hence, we conclude that $\min(A_9(C_{59,1}))=5$,
since $\mu(59)=5$.
For $A_9(C_{59,1})$, we obtained $(\beta,\gamma,\tau)=(0,6844,17228)$.

We also found $10$ modified double circulant self-dual $\ZZ_9$-codes
$C_{59,i}$ of length $59$ such that $A_9(C_{59,i})$ has
minimum norm $5$ for $i=2,3,\ldots,11$.
The codes $C_{59,i}$ are generated by the rows of matrices of the
form~\eqref{eq:odd-matrix}, where $(a,b,c)=(2,0,3)$ and
the first rows $r$ of the $29\times29$ circulant matrices $R$ are
listed in Table~\ref{Tab:59Z9}.

By Proposition~\ref{prop:rank59-min5-theta}, the theta series of
each $A_9(C_{59,i})$ is determined by the parameters $\beta$ and
$\gamma$ defined there. 
For each $i$ with $2\leq i\leq11$, Table~\ref{Tab:59Z9} lists
the first row $r$ of $R$ and the triple $(\beta,\gamma,\tau)$.
For the lattice in dimension $59$ constructed in~\cite[Table~1]{SYS},
we obtained $(\beta,\gamma,\tau)=(0,6136,20060)$.
The $11$ lattices $A_9(C_{59,i})$
and the lattice constructed in~\cite{SYS} have distinct
kissing numbers. Since $\mu(59)=5$, we obtain the following:

\begin{prop}\label{prop:L59}
There exist at least $12$ nonisomorphic optimal unimodular lattices
with minimum norm $5$ in dimension $59$.
\end{prop}

\begin{table}[thb]
\caption{Unimodular lattices $A_9(C_{59,i})$ $(i=2,3,\ldots,11)$}
\label{Tab:59Z9}
\centering
\medskip
{\footnotesize
\begin{tabular}{c|c|c}
\noalign{\hrule height1pt}
$i$ & $r$ & $(\beta,\gamma,\tau)$ \\
\hline
2&$(7,3,1,1,2,2,2,3,2,2,2,1,2,2,1,2,3,1,2,1,3,1,2,3,3,3,2,2,2)$&$(0,7076,16300)$\\
3&$(7,7,3,7,5,7,0,7,0,1,6,4,7,8,4,5,4,8,2,4,5,0,0,4,7,4,8,2,0)$&$(0,5916,20940)$\\
4&$(8,7,5,0,0,7,6,4,2,3,2,2,0,2,2,6,2,1,1,5,6,7,1,0,1,8,2,0,0)$&$(0,5684,21868)$\\
5&$(6,6,6,0,1,7,2,3,2,5,4,1,6,2,0,5,6,6,6,5,8,1,2,3,0,6,2,2,5)$&$(0,5452,22796)$\\
6&$(0,8,0,2,3,0,8,1,3,1,3,0,6,6,6,8,5,8,6,4,1,2,2,5,8,8,0,4,0)$&$(0,5336,23260)$\\
7&$(8,8,2,2,2,6,1,8,0,6,8,8,6,0,1,4,6,4,8,4,3,8,6,6,2,7,6,7,7)$&$(0,5220,23724)$\\
8&$(4,5,3,1,5,4,1,7,5,6,8,8,7,0,2,3,5,6,1,7,3,6,1,4,1,4,5,1,4)$&$(0,5104,24188)$\\
9&$(0,8,7,2,8,8,5,2,2,0,8,2,7,4,2,5,4,2,2,3,2,1,0,4,3,0,7,0,1)$&$(0,4872,25116)$\\
10&$(1,8,7,6,7,2,6,1,0,0,6,5,1,0,4,0,8,4,1,6,0,7,2,5,8,7,0,5,1)$&$(0,4756,25580)$\\
11&$(8,0,5,5,4,4,4,0,4,4,4,5,4,4,5,4,0,5,4,5,0,5,4,0,0,0,4,4,4)$&$(0,4524,26508)$\\
\noalign{\hrule height1pt}
\end{tabular}
}
\end{table}

\subsection{Table of $\mu(n)$}

As noted in Section~\ref{Sec:1}, $\mu(n)$ is known for every dimension $n\leq48$.
Combining the results of this paper with those in~\cite{Ga04},
\cite{NG} and~\cite{SYS}, we obtain the values $\mu(n)$
for $49\leq n\leq61$ listed in Table~\ref{Tab:mu}.
Throughout this range, $\mu(n)=5$ except in dimensions $49$, $56$
and $60$, where $\mu(49)=4$ and $\mu(56)=\mu(60)=6$.
Therefore, $\mu(n)$ is now known for every dimension $n\leq61$.

\begin{table}[thb]
\caption{Largest minimum norms  $\mu(n)$ $(n=49,50,\ldots,61)$}
\label{Tab:mu}
\centering
\medskip
{\small
\begin{tabular}{c|c||c|c||c|c}
\noalign{\hrule height1pt}
$n$ & $\mu(n)$ &$n$ & $\mu(n)$ & $n$ & $\mu(n)$ \\
\hline
$49$ & $4$      & $54$ & $5$      &   $59$ & $5$    \\
$50$ & $5$      & $55$ & $5$      &  $60$ & $6$    \\
$51$ & $5$      & $56$ & $6$      &   $61$ & $5$     \\
$52$ & $5$      & $57$ & $5$      &      \\
$53$ & $5$      &  $58$ & $5$          &          \\
\noalign{\hrule height1pt}
\end{tabular}
}
\end{table}

\bigskip
\noindent
\textbf{Acknowledgments.}
This work was supported by JSPS KAKENHI Grant Numbers 23K25784 and
25K21995.
The author would like to thank Patrick Sol\'e for bringing the recent preprint \cite{SYS}
to the author's attention.

\bigskip
\noindent
\textbf{Declaration of Generative AI Use.}
In November 2026,
the author used OpenAI's ChatGPT (GPT-6 Astra) to assist in developing the proofs of 
Proposition~\ref{prop:rank49-min4-theta} and
Theorem~\ref{thm:common-upper-bounds}.
This assistance also extended to the proofs of three supporting lemmas
(Lemmas~\ref{lem:theta-common-coefficients},
\ref{lem:theta-shadow-bounds} and~\ref{lem:theta-weighted-identities}).
The author independently verified all mathematical arguments
presented in the paper.



\end{document}